\documentclass[12pt]{article}

\usepackage{vmargin}
\setpapersize{A4}
\setmargins{2.5cm}       
{1.5cm}                        
{16.5cm}                      
{23.42cm}                    
{10pt}                           
{1cm}                           
{0pt}                             
{2cm}                           
\usepackage{ifpdf}
\usepackage[mathscr]{eucal}
\usepackage{amsfonts}
\usepackage{amssymb}
\usepackage{color}
\usepackage {graphics}
\usepackage {epsfig}
        \usepackage{latexsym}
        \usepackage{palatino}
        \usepackage{fancybox}
\usepackage{amsmath,mathptmx,amssymb,bm}
\usepackage{amsthm}

\newtheorem{example}{Example}

\newtheorem{rem}{Remark}

\usepackage{titlesec}

\titleformat{\paragraph}
{\normalfont\normalsize\bfseries}{\theparagraph}{1em}{}
\titlespacing*{\paragraph}
{0pt}{3.25ex plus 1ex minus .2ex}{1.5ex plus .2ex}

\begin{document}

\title{\textbf{The natural extension to PDEs of Lie's reduction of order algorithm for ODEs}}

\author{ George W. Bluman${}^{a)}$, Rafael de la Rosa${}^{b)}$ \\ 
 ${}^{a)}$ Department of Mathematics, University of British Columbia, Vancouver, \\ British Columbia V6T 1Z2, Canada (e-mail:bluman@math.ubc.ca)\\
${}^{b)}$ Departamento de Matem\'aticas, Universidad de C\'adiz, 11510 Puerto Real,\\ C\'adiz, Spain (e-mail:rafael.delarosa@uca.es)}

\date{}
 
\maketitle

\begin{abstract}

In this paper, we further consider the symmetry-based method for seeking nonlocally related systems for partial differential equations. In particular, we show that the symmetry-based method for partial differential equations is the natural extension of Lie's reduction of order algorithm for ordinary differential equations by looking at this algorithm from a different point of view. Many examples exhibit various situations that can arise.


\end{abstract}






\section{Introduction}

A differential equation (DE) is a ``compact equation'' for describing a family of curves (surfaces) that are the solutions of the DE. The curves (surfaces) are continuously connected. A symmetry of a DE maps any solution of the DE to another solution of the same DE. Hence, from a thought experiment point of view, every DE has symmetries.

One can seek symmetries of a DE without knowing its solutions.  In turn, for a given DE, one can use symmetries to find specific solutions, seek conservation laws if the DE is variational, or seek a mapping of the DE to a simpler DE (e.g., mapping a linear partial differential equation (PDE) with variable coefficients to a linear PDE with constant coefficients, mapping a nonlinear PDE to a linear PDE).

In practice, a useful symmetry needs to be described in some coordinate space, not necessarily the coordinate space of a considered DE. A Lie point symmetry is a useful symmetry that acts on the space of independent and dependent variables of a given DE.

Sophus Lie introduced the notion of a one-parameter group of point transformations, now known as a one-parameter Lie group of point transformations, in order to systematically unify and extend specialized, mostly ad-hoc, methods for solving ordinary differential equations (ODEs). Lie presented an algorithm (Lie's algorithm) to seek one-parameter Lie groups of point transformations that leave invariant a given DE (an ODE or PDE) without knowing the solutions of the DE. Such symmetries are called Lie point symmetries or, more simply, point symmetries, of a DE.  Lie's algorithm showed that the order of a given ODE can be constructively reduced by one when the given ODE is invariant under a point symmetry. Furthermore, Lie showed that solutions of the given ODE can be recovered from solutions of the reduced ODE through quadrature. For PDEs, Lie \cite{lie1881} showed that a point symmetry yields special solutions called invariant or similarity solutions. Lie's fundamental theorems showed that point symmetries are completely characterized by their infinitesimal generators.

The applicability of symmetry methods can be further enhanced through systematically seeking nonlocal symmetries. Here one  seeks an equivalent DE system, which is not invertibly related to the given DE system. Such an equivalent DE system is called a nonlocally related system. Each solution of the given DE system yields a solution of the nonlocally related
system and, vice versa, any solution of the nonlocally related system yields a solution of the given DE system. But there is not a one-to-one relationship between the solutions of these two equivalent systems. Furthermore, there needs to be a connection formula that relates the solutions of the nonlocally related systems. More importantly, useful nonlocal symmetries for a given DE system could arise from seeking point symmetries of a nonlocally related system by applying Lie's algorithm. Thus, nonlocally related systems can play an essential role in the analysis of a given DE system.

There are two systematic approaches for the construction of nonlocally related systems: the conservation law (CL)-based method and the symmetry-based method. For both methods, there is a connection formula that relates all solutions of a nonlocally related DE system with all solutions of the given DE system.

The CL-based method was introduced in \cite{blumankumei1987} and further developed in \cite{blumankumeireid1988}. In \cite{blumankumeireid1988}, a systematic method was presented for seeking nonlocally related systems and consequent nonlocal symmetries (potential symmetries) from a CL of a given PDE system. This method was extended to include PDE systems admitting more than one CL and PDE systems with more than two independent variables. See references \cite{blumankumei,ancobluman1997,ancothe2005,bluman2006,bluman2007,blumanchevianco} for further details.

The symmetry-based method was introduced in \cite{blumanyang2013,yang2013}. In this method, a nonlocally
related PDE system is systematically constructed from any admitted point symmetry of a given PDE system. Here, taking into account the canonical coordinates of an admitted point symmetry with the translated canonical coordinate playing the role of a new dependent variable, one can construct a nonlocally related PDE system. The symmetry-based method is especially useful when a given PDE system has no CLs but does admit point symmetries.

In \cite{DIpaper2}, we proved two important results:
\begin{itemize}
\item[(1)] For a PDE system with two independent variables, the symmetry-based method includes the CL-based method as a special case.
\item[(2)] For a PDE system with three or more independent variables, the symmetry-based method includes the
CL-based method for curl-type CLs as a special case.
\end{itemize}

In \cite{DIpaper1} we showed that there is a connection between Lie's reduction of order method for ODEs and the symmetry-based method for seeking nonlocally related systems for PDEs. In the current paper, we further explore this connection. By looking at Lie's reduction of order algorithm for an ODE from a different point of view, we see that a given ODE is nonlocally related to its reduced ODE.  More importantly, we show that the symmetry-based method for finding nonlocally related PDE systems of a PDE with any number of independent variables is the natural extension to PDEs of Lie's reduction of order algorithm for ODEs.  

The rest of the paper is organized as follows.

In Section \ref{sec:liereducODE}, we consider Lie's reduction of order algorithm from a different point of view to show how the given ODE is nonlocally related to its reduced ODE. We give examples illustrating how the relationship between the solutions of these equivalent systems is not one-to-one. Three different symmetry situations arise that are illustrated in the examples: a second point symmetry of a given ODE is also a point symmetry of the reduced ODE; a second point symmetry of a given ODE is a symmetry that is not a point symmetry of the reduced ODE; a point symmetry of the reduced ODE is a symmetry that is not a point symmetry of the given ODE.

In Section \ref{sec:PDE}, we show that the natural extension to PDEs of Lie's reduction of order method for ODEs is simply the symmetry-based method for obtaining a nonlocally related system for a PDE with any number of independent variables. We consider separately the situation for PDEs with two independent variables, three independent variables and more than three independent variables. Analogously, we give examples illustrating how the relationship between the solutions of these equivalent systems is not one-to-one, as well as illustrate the same three different symmetry situations that can arise.

Finally, in Section \ref{sec:concluding}, we present some important concluding remarks. These include the extension of the theory presented in Sections \ref{sec:liereducODE} and \ref{sec:PDE} to systems of DEs with more than one dependent variable; an explanation of what ``reduction of order'' means for PDE systems admitting a point symmetry; the natural extension of the use of solvable Lie algebras in the reduction of order algorithm for ODEs to the use of solvable Lie algebras in the symmetry-based method for PDEs; the importance of nonlocally related systems for enhancing the use of various methods for DEs.


\section{Lie's reduction of order algorithm for an ODE}\label{sec:liereducODE}

Consider an $n$th-order ODE with independent variable $x$ and dependent variable $y$ given by
\begin{equation}
\label{ODE}
K\left(x,y,y',\ldots,y^{(n)}\right)=0.
\end{equation}
Assume that ODE (\ref{ODE}) admits a point symmetry with infinitesimal generator
\begin{equation}
\label{symmODE}
\mathbb{X}=\xi\left(x,y\right)\dfrac{\partial}{\partial x}+\eta\left(x,y\right)\dfrac{\partial}{\partial y}
\end{equation}
whose $n$th extension is given by
\begin{equation}
\label{symmODEext}
\mathbb{X}^{(n)}=\mathbb{X}+ \eta^{(1)}\left(x,y,y'\right)\dfrac{\partial}{\partial y'}+\cdots+\eta^{(n)}\left(x,y,y',\ldots,y^{(n)}\right)\dfrac{\partial}{\partial y^{(n)}}
\end{equation}
with $\eta^{(k)}=D \eta^{(k-1)}-y^{(k)}D \xi$, for $k=1,\ldots,n$, $\eta^{(0)}=\eta(x,y)$, in terms of the total differential operator
$$D=\frac{\partial}{\partial x}+y' \frac{\partial}{\partial y}+\cdots+y^{(k)}\frac{\partial}{\partial y^{(k-1)}}+\ldots \, \,.$$
Then $\mathbb{X}$ is a point symmetry of ODE (\ref{ODE}) if and only if
$$ \mathbb{X}^{(n)}K\left(x,y,y',\ldots,y^{(n)}\right) =0 \quad \mbox{when} \quad K\left(x,y,y',\ldots,y^{(n)}\right)=0.$$
We introduce canonical coordinates, $r(x,y)$ and $s(x,y)$, for $\mathbb{X}$, satisfying
$$\mathbb{X}r=0, \quad \mathbb{X}s=1.$$
Thus, equivalently, ODE (\ref{ODE}) is invariant under translations in $s$, i.e., ODE (\ref{ODE}) is invariant under the point symmetry (\ref{symmODE}) now written as $\mathbb{X}=\dfrac{\partial}{\partial s}$.\vspace*{0.1cm}

In terms of corresponding canonical coordinates, ODE (\ref{ODE}) invertibly transforms by a point transformation to an ODE of the form
\begin{equation}
\label{ODE2}
\widehat{K}\left(r,\dfrac{ds}{dr},\ldots,\dfrac{d^n s}{dr^n}\right)=0.
\end{equation}

Introducing the auxiliary dependent variable $z=\dfrac{ds}{dr}$, ODE (\ref{ODE}) reduces to the $(n-1)$th-order ODE
\begin{equation}
\label{reducODE}
\widehat{K}\left(r,z,\dfrac{dz}{dr},\ldots,\dfrac{d^{n-1} z}{dr^{n-1}}\right)=0.
\end{equation}

If $z=Z(r)$ is a solution of ODE (\ref{reducODE}), then
$$s=s(r(x,y))=\int^{r(x,y)} Z(\rho) \, d\rho+C$$
solves ODE (\ref{ODE}) for any constant $C$. For more detailed information on Lie's reduction of order algorithm see any of the books \cite{page1897,cohen1911,friedrichs1965,blumancole1974,ovsiannikov1982,olver1986,blumankumei}.

Now we look at Lie's reduction of order algorithm for reducing the order of an ODE from a different point of view. Without loss of generality, if a given ODE admits a point symmetry, after using corresponding canonical coordinates, it can be re-written in the form
\begin{equation}
\label{newODE}
L\left(x,y',y'',\ldots,y^{(n)}\right)=0,
\end{equation}
i.e., ODE (\ref{newODE}) is invariant under the translation symmetry $\mathbb{X}=\dfrac{\partial}{\partial y}$.\vspace*{0.1cm}

From the form of ODE (\ref{newODE}), we can naturally introduce as an auxiliary dependent variable $\alpha=y'$. Then ODE (\ref{newODE}) becomes the reduction of order ODE
\begin{equation}
\label{reducnewODE}
L\left(x,\alpha,\alpha',\ldots,\alpha^{(n-1)}\right)=0, 
\end{equation}
nonlocally related to ODE (\ref{newODE}).\\

\begin{rem}\label{rem:solsODE}
\noindent If 
$\alpha=A(x)$ solves ODE (\ref{reducnewODE}), then
$$y=\int^{x} A(\rho)\, d\rho+C$$
solves ODE (\ref{newODE}) for any constant $C$. Conversely, if $y=Y(x)$ solves ODE (\ref{newODE}), then $\alpha=Y'(x)$ solves ODE (\ref{reducnewODE}). But the relationship between the solutions of ODEs (\ref{newODE}) and (\ref{reducnewODE}) is not one-to-one.

In particular, for any solution of ODE (\ref{reducnewODE}), there is an infinite number of solutions of ODE (\ref{newODE}). Conversely, for any solution of ODE (\ref{newODE}), there is a unique solution of ODE (\ref{reducnewODE}). Thus ODEs (\ref{newODE}) and (\ref{reducnewODE}) are nonlocally related.\\
\end{rem}

We consider an example that appeared in \cite{olver1986}.\\

\begin{example}
\normalfont
Consider the second-order ODE given by
\begin{equation}\label{canonolverODE}
y''=(1+x) \left(y'\right)^2+y'.
\end{equation}
Let $\alpha=y'$. Then ODE (\ref{canonolverODE}) reduces to the nonlocally related first-order ODE (Bernoulli equation) 
\begin{equation}\label{reduccanonolverODE}
\dfrac{d \alpha}{d x}=(1+x) \alpha^2+\alpha.
\end{equation}
Consider a particular solution of ODE (\ref{reduccanonolverODE}) given by
\begin{equation}\label{solbernoulliODE}
\alpha=\dfrac{1}{e^{-x}-x}.
\end{equation}
Correspondingly
\begin{equation}\label{sololverODE}
y=\int^{x}\dfrac{1}{
e^{-\rho}-\rho}\, d\rho+C
\end{equation}
solves ODE (\ref{canonolverODE}) with $C$ an arbitrary constant. 

Now consider a particular solution of ODE (\ref{canonolverODE}) given by
\begin{equation}\label{solcanonolverODE}
y=-\log x.
\end{equation}
Correspondingly
\begin{equation}\label{sol2bernoulliODE}
\alpha=-\dfrac{1}{x}
\end{equation}
solves ODE (\ref{reduccanonolverODE}). This illustrates that the relationship between the solutions of ODE (\ref{canonolverODE}) and ODE (\ref{reduccanonolverODE}) is not one-to-one.\\
\end{example}

\begin{rem}\label{rem:symmODE}
\noindent
From the relationship between the solutions of ODE (\ref{newODE}) and ODE (\ref{reducnewODE}), any symmetry of ODE (\ref{newODE}) is a symmetry of ODE (\ref{reducnewODE}); and vice versa, any symmetry of ODE (\ref{reducnewODE}) is a symmetry of ODE (\ref{newODE}).

A point symmetry of ODE (\ref{newODE}) is not necessarily a point symmetry of ODE (\ref{reducnewODE}). Conversely, a point symmetry of ODE (\ref{reducnewODE}) is not necessarily a point symmetry of ODE (\ref{newODE}).\\
\end{rem}

\begin{example}\label{ex:nonlocalsymmgivenODE}
\normalfont
Consider again ODE (\ref{canonolverODE}). It can be easily checked the only point symmetry that ODE (\ref{canonolverODE}) admits is the translation symmetry $\mathbb{X}=\frac{\partial}{\partial y}$. Since ODE (\ref{reduccanonolverODE}) is a first-order ODE, it admits an infinite number of point symmetries. In particular, ODE (\ref{reduccanonolverODE}) admits the point symmetry generator 
\begin{equation}\label{newsymmreduccanonolverODE}
\mathbb{Y}=\alpha (1+x \alpha) \dfrac{\partial}{\partial \alpha}.
\end{equation}

Clearly, point symmetry (\ref{newsymmreduccanonolverODE}) is a symmetry of ODE (\ref{canonolverODE}) that is not a point symmetry of ODE (\ref{canonolverODE}) since ODE (\ref{canonolverODE}) only admits one point symmetry.\\
\end{example}
\medskip


To illustrate the other two of the three different symmetry situations that could arise, we consider a second-order ODE admitting two point symmetries such that: (1) the second point symmetry yields a point symmetry of the reduced ODE obtained from its invariance under the first point symmetry; (2) the second point symmetry yields a symmetry which is not a point symmetry of the reduced ODE obtained from its invariance under the first point symmetry (this symmetry is a nonlocal symmetry of the reduced ODE). [The third symmetry situation that could arise has been illustrated in Example \ref{ex:nonlocalsymmgivenODE}: (3) a point symmetry of the reduced ODE obtained from the invariance under a point symmetry of the given ODE yields a symmetry which is not a point symmetry of the given ODE (this symmetry is a nonlocal symmetry of the given ODE).]\\

\begin{example} \label{ex:ODEbothdirec}
\normalfont
Consider the second-order ODE given by
\begin{equation}\label{secondex}
x y^2 y''+x y'-y=0.
\end{equation}
ODE (\ref{secondex}) admits two point symmetries given by
\begin{equation}\label{symmsecondex} \mathbb{X}_1=x^2\dfrac{\partial}{\partial x}+x y\dfrac{\partial}{\partial y}, \quad \mathbb{X}_2=x\dfrac{\partial}{\partial x}+\dfrac{1}{2}y\dfrac{\partial}{\partial y}.
\end{equation}
\medskip

\noindent Consider $\mathbb{X}_1$. In terms of canonical coordinates
\begin{equation}
\label{canon1secondex} r=\dfrac{y}{x}, \qquad s=-\dfrac{1}{x},
\end{equation}
ODE (\ref{secondex}) becomes
\begin{equation}\label{canon1secondexODE}
r^2 \dfrac{d^2 s}{d r^2}-\left(\dfrac{ds}{dr}\right)^2=0.
\end{equation}
Let
\begin{equation}\label{auxvar1secondex} \alpha=\dfrac{ds}{dr}=\dfrac{1}{x y'-y}.
\end{equation}
Then ODE (\ref{secondex}) reduces to the nonlocally related separable first-order ODE 
\begin{equation}\label{reduccanon1secondexODE}
r^2 \dfrac{d \alpha}{d r}-\alpha^2=0.
\end{equation}
After applying 
$$\mathbb{X}_2^{(1)}=\mathbb{X}_2-\dfrac{1}{2}y' \dfrac{\partial}{\partial y'},$$
to (\ref{canon1secondex}) and (\ref{auxvar1secondex}), one obtains
$$\mathbb{X}_2^{(1)}r=-\dfrac{1}{2}r, \qquad \mathbb{X}_2^{(1)}\alpha=-\dfrac{1}{2}\alpha.$$
Hence the point symmetry $\mathbb{X}_2$ in (\ref{symmsecondex}) of ODE (\ref{secondex}) becomes
\begin{equation}\label{newX2}\mathbb{X}_2=r \dfrac{\partial}{\partial r}+\alpha \dfrac{\partial}{\partial \alpha}.
\end{equation}
The symmetry (\ref{newX2}) of ODE (\ref{reduccanon1secondexODE}) is a point symmetry of ODE (\ref{reduccanon1secondexODE}) since symmetry (\ref{newX2}) acts on $\left(r,\alpha\right)$-space.\\

\noindent Now consider $\mathbb{X}_2$. In terms of canonical coordinates
\begin{equation}
\label{canon2secondex} r=\dfrac{y^2}{x}, \qquad s=\log x,
\end{equation}
ODE (\ref{secondex}) becomes
\begin{equation}\label{canon2secondexODE}
2 r \dfrac{d^2 s}{d r^2}+ r (r+2)\left(\dfrac{ds}{dr}\right)^3-2 \left(\dfrac{ds}{dr}\right)^2+\dfrac{ds}{dr}=0.
\end{equation}
Let
\begin{equation}\label{auxvar2secondex} \alpha=\dfrac{ds}{dr}=\dfrac{x}{2 x y y'-y^2}.
\end{equation}
Then ODE (\ref{secondex}) reduces to the nonlocally related first-order ODE (Abel equation of the first kind) 
\begin{equation}\label{reduccanon2secondexODE}
2 \alpha \dfrac{d \alpha}{d r}+ r (r+2)\alpha^3-2 \alpha^2+\alpha=0.
\end{equation}
After applying 
$$\mathbb{X}_1^{(1)}=\mathbb{X}_1+\left(y-x y'\right)\dfrac{\partial}{\partial y'},$$
to (\ref{canon2secondex}) and (\ref{auxvar2secondex}), one obtains
$$\mathbb{X}_1^{(1)}r=r e^{s}, \qquad \mathbb{X}_1^{(1)}\alpha=-r e^{s} \alpha^2.$$
Hence the point symmetry $\mathbb{X}_1$ in (\ref{symmsecondex}) of ODE (\ref{secondex}) becomes the symmetry
\begin{equation}\label{newX1}\mathbb{X}_1=r e^{s}\dfrac{\partial}{\partial r}-r e^{s} \alpha^2 \dfrac{\partial}{\partial \alpha},
\end{equation}
of ODE (\ref{reduccanon2secondexODE}).

But although ODE (\ref{reduccanon2secondexODE}) is a first-order ODE and thus admits an infinite number of point symmetries, the symmetry (\ref{newX1}) of ODE (\ref{reduccanon2secondexODE}) is not a point symmetry of ODE (\ref{reduccanon2secondexODE}) since $\mathbb{X}_1$ depends explicitly on $s$. In particular $\mathbb{X}_1$ is a nonlocal symmetry of ODE (\ref{reduccanon2secondexODE}).\\
\end{example}

To further illustrate Remark \ref{rem:symmODE}, we consider the relationship between the number of point symmetries admitted by an $n$th-order ODE and the reduced $(n-1)$th-order ODE obtained from invariance under a point symmetry of the $n$th-order ODE. Note that (see \cite{blumananco,dickson1924,lie1893,ovsiannikov1982} for further details):
\begin{enumerate}
\item Any first-order ODE admits an infinite number of nontrivial point symmetries.
\item A second-order ODE admits at most eight point symmetries.
\item The number of point symmetries admitted by an $n$th-order ODE ($n \geq 3$) is at most $n + 4$.
\end{enumerate}
Hence, the reduced first-order ODE arising from a second-order ODE always has more point symmetries than the original ODE. The reduced ODE arising from a third-order ODE could have more point symmetries than the original ODE. On the other hand, the reduced ODE arising from an $n$th-order ODE ($n \geq 4$) could have fewer point symmetries than the original ODE.\\

\begin{rem}\label{rem:furtherNLRODE}
\noindent
If a point symmetry of ODE (\ref{newODE}) is also a point symmetry of ODE (\ref{reducnewODE}), one can construct explicitly another nonlocally related ODE for ODE (\ref{reducnewODE}) in addition to ODE (\ref{newODE}), and a further nonlocally related ODE for ODE (\ref{newODE}).\\
\end{rem}

\begin{example} \label{ex:ODEbothdireccont}
\normalfont
Consider again the second-order ODE (\ref{secondex}) and its nonlocally related reduced first-order ODE (\ref{reduccanon1secondexODE}). Consider $\mathbb{X}_2$ given by (\ref{newX2}). From Example \ref{ex:ODEbothdirec} we know that $\mathbb{X}_2$ is a point symmetry of ODE (\ref{reduccanon1secondexODE}). In terms of canonical coordinates
\begin{equation}
\label{furthercanon1secondex} R=\dfrac{\alpha}{r}, \qquad S=\log r,
\end{equation}
ODE (\ref{reduccanon1secondexODE}) becomes
\begin{equation}\label{furthercanon1secondexODE}
1+R(1-R)\dfrac{d S}{d R}=0.
\end{equation}
Let $\omega=\dfrac{dS}{dR}$. Then ODE (\ref{reduccanon1secondexODE}) reduces to the algebraic equation
\begin{equation}\label{reducfurthercanon1secondexODE}
1+R(1-R)\omega=0,
\end{equation}
nonlocally related to ODE (\ref{furthercanon1secondexODE}) and ODE (\ref{secondex}).\\
\end{example}

\begin{example}\label{ex:blasius}
\normalfont
As a second example, consider the Blasius equation
\begin{equation}
\label{blasius} y'''+\frac{1}{2} y y''=0.
\end{equation} 
ODE (\ref{blasius}) admits two point symmetries given by
\begin{equation}
\label{symmblasius} \mathbb{X}_1=\dfrac{\partial}{\partial x}, \quad \mathbb{X}_2=x\dfrac{\partial}{\partial x}-y\dfrac{\partial}{\partial y}.
\end{equation}
\medskip

\noindent Consider $\mathbb{X}_1$. In terms of canonical coordinates
\begin{equation}
\label{canonvarblasius} r=y, \qquad s=x,
\end{equation}
the Blasius equation (\ref{blasius}) becomes the invertibly related ODE
\begin{equation}\label{canonblasius}
2 \dfrac{ds}{dr} \dfrac{d^3 s}{d r^3}-6 \left( \dfrac{d^2 s}{d r^2}\right)^2+ r  \left(\dfrac{d s}{d r}\right)^2 \dfrac{d^2 s}{d r^2}=0.
\end{equation}
Without loss of generality, ODE (\ref{canonblasius}) can be re-written in the form
\begin{equation}\label{canonblasiusorig}
2 y' y'''-6 \left( y'' \right)^2+ x  \left(y'\right)^2 y''=0.
\end{equation}
ODE (\ref{canonblasiusorig}) obviously admits the point symmetries given by 
\begin{equation}
\label{symmblasiusorig} \mathbb{X}_1=\dfrac{\partial}{\partial y}, \quad \mathbb{X}_2=x\dfrac{\partial}{\partial x}-y\dfrac{\partial}{\partial y}.
\end{equation}
Let
\begin{equation}\label{auxvar1blasius} \alpha=y'.
\end{equation}
Then ODE (\ref{canonblasiusorig}) reduces to the second-order ODE 
\begin{equation}\label{reduccanonblasius}
2 \alpha \dfrac{d^2 \alpha}{d x^2}-6 \left( \dfrac{d \alpha}{d x}\right)^2+ x  \alpha^2 \dfrac{d \alpha}{d x}=0,
\end{equation}
nonlocally related to ODE (\ref{canonblasiusorig}).

Then
$$\mathbb{X}_2^{(1)}=\mathbb{X}_2-2y' \dfrac{\partial}{\partial y'}.$$
Correspondingly, the point symmetry $\mathbb{X}_2$ in (\ref{symmblasiusorig}) of ODE (\ref{canonblasiusorig}) becomes the point symmetry 
\begin{equation}\label{newX2blasius}\mathbb{X}_2=x \dfrac{\partial}{\partial x}-2\alpha \dfrac{\partial}{\partial \alpha},
\end{equation}
of ODE (\ref{reduccanonblasius}).\\

\noindent Now consider $\mathbb{X}_2$. In terms of canonical coordinates
\begin{equation}
\label{furthercanon1blasius} r=x^2 \alpha, \qquad s=\log x,
\end{equation}
ODE (\ref{reduccanonblasius}) becomes
\begin{equation}\label{furthercanon1blasiusODE}
2 r \dfrac{d^2 s}{d r^2}+2 r^2 (r+6) \left(\dfrac{d s}{d r} \right)^3-r (r+14)\left(\dfrac{d s}{d r} \right)^2+6 \dfrac{d s}{d r}=0.
\end{equation}
Let $\omega=\dfrac{ds}{dr}$. Then ODE (\ref{reduccanonblasius}) reduces to the first-order ODE
\begin{equation}\label{reducfurthercanon1blasiusODE}
2 r \dfrac{d \omega}{d r}+2 r^2 (r+6) \omega^3-r (r+14)\omega^2+6 \omega=0,
\end{equation}
nonlocally related to ODE (\ref{reduccanonblasius}) and ODE (\ref{canonblasiusorig}).
\end{example}

\section{The natural extension to a PDE of Lie's reduction of order algorithm for an ODE}\label{sec:PDE}

In this section, we present the natural extension of Lie's reduction of order algorithm for an ODE to a PDE with any number of independent variables. We show how to do this through further consideration of the symmetry-based method to obtain nonlocally related systems/nonlocal symmetries of a given PDE system \cite{blumanyang2013,yang2013,blumanyang2014,blumanzuhal,DIpaper1,DIpaper2,wang2023}. We consider separately the cases of PDEs with two independent variables, three independent variables and more than three independent variables.

\subsection{A PDE with two independent variables}\label{pde2}
Consider an $n$th-order PDE with two independent variables $x_1$ and $x_2$ and dependent variable $u$ given by
\begin{equation}
\label{PDE}
K\left(x_1,x_2,u,\partial u,\ldots,\partial^n u \right)=0,
\end{equation}
with $\partial^l u$ denoting its partial derivatives of $u$ of order $l$, $l=1,2,\ldots,n$.

Assume that PDE (\ref{PDE}) admits a point symmetry with infinitesimal generator
\begin{equation}
\label{symmPDE}
\mathbb{X}=\xi_1 \left(x_1,x_2,u\right)\dfrac{\partial}{\partial x_1}+\xi_2 \left(x_1,x_2,u\right)\dfrac{\partial}{\partial x_2}+\eta\left(x_1,x_2,u\right)\dfrac{\partial}{\partial u}
\end{equation}
whose $n$th extension is given by
\begin{equation}
\label{symmPDEext}
\mathbb{X}^{(n)}=\mathbb{X}+ \eta_j^{(1)}\left(x_1,x_2,u,\partial u\right)\dfrac{\partial}{\partial u_j}+\cdots+\eta_{j_1 j_2 \cdots j_n}^{(n)}\left(x_1,x_2,u,\partial u,\ldots,\partial^n u \right)\dfrac{\partial}{\partial u_{j_1 j_2 \cdots j_n}},
\end{equation}
where $u_j=\frac{\partial u}{\partial x_j}$, $j=1,2$, and $j_k=1,2,$ for $k=1,\ldots,n$, and the $k$th-extended infinitesimals  $\eta^{(k)}$ satisfy the recursion relations
$$\begin{array}{rcl}\eta_j^{(1)}&=& D_j \eta- \left( D_j \xi_i \right)u_i, \quad j=1,2, \vspace*{0.05cm}\\ 
\eta_{j_1 j_2 \cdots j_k}^{(k)}&=&  D_{j_k} \eta_{j_1 j_2 \cdots j_{k-1}}^{(k-1)}- \left( D_{j_k} \xi_i \right)u_{j_1 j_2 \cdots j_{k-1} i},
\end{array}$$
in terms of the total derivative operators
\begin{equation}\label{totalderivative} D_j=\dfrac{\partial}{\partial x_j}+u_j \dfrac{\partial}{\partial u}+u_{j j_1} \dfrac{\partial}{\partial u_{j_1}}+\ldots+u_{j j_1 \ldots j_k} \dfrac{\partial}{\partial u_{j_1 \ldots j_k}}+ \ldots \, \, .\end{equation}
Then $\mathbb{X}$ is a point symmetry of PDE (\ref{PDE}) if and only if 
$$ \mathbb{X}^{(n)}K\left(x_1,x_2,u,\partial u,\ldots,\partial^n u \right) =0 \quad \mbox{when} \quad K\left(x_1,x_2,u,\partial u,\ldots,\partial^n u \right)=0.$$

For details on how to find and use point symmetries for PDEs see any of the books \cite{ovsiannikov1962,ovsiannikov1982,olver1986,
blumankumei,stephani1989,hydon2000,blumananco,cantwell2002,blumanchevianco}.

We introduce canonical coordinates, $r_1(x_1,x_2,u)$, $r_2(x_1,x_2,u)$ and $s(x_1,x_2,u)$, for $\mathbb{X}$, satisfying
$$\mathbb{X}r_1=0, \quad \mathbb{X}r_2=0, \quad  \mathbb{X}s=1.$$
Thus, equivalently, PDE (\ref{PDE}) is invariant under translations in $s$, i.e., PDE (\ref{PDE}) is invariant under the point symmetry (\ref{symmPDE}) now written as $\mathbb{X}=\dfrac{\partial}{\partial s}$.\vspace*{0.1cm}

In terms of corresponding canonical coordinates, PDE (\ref{PDE}) invertibly transforms by a point transformation to a PDE of the form
\begin{equation}
\label{transformPDE}
\widehat{K}\left(r_1,r_2,\partial s,\ldots,\partial^n s \right)=0.
\end{equation}

Hence, without loss of generality, if a given PDE with two independent variables admits a point symmetry, after using corresponding canonical coordinates, it can be re-written in the form
\begin{equation}
\label{newPDE}
L\left(x_1,x_2,\partial u,\ldots,\partial^n u \right)=0,
\end{equation}
i.e., PDE (\ref{newPDE}) is invariant under the translation symmetry $\mathbb{X}=\dfrac{\partial}{\partial u}$.\vspace*{0.1cm}

From the form of PDE (\ref{newPDE}), we can naturally introduce as auxiliary dependent variables $\alpha=\frac{\partial u}{\partial x_1}$, $\beta=\frac{\partial u}{\partial x_2}$. Then PDE (\ref{newPDE}) becomes the ``reduction of order'' PDE system
\begin{equation}
\label{reducnewPDE}
\left\{ \begin{array}{l} L\left(x_1,x_2,\alpha, \beta,\partial \alpha, \partial \beta, \ldots,\partial^{n-1} \alpha, \partial^{n-1} \beta \right) =0, \\ \vspace*{0.15cm}
\frac{\partial \alpha}{\partial x_2}=\frac{\partial \beta}{\partial x_1},
\end{array} \right.
\end{equation}
nonlocally related to PDE (\ref{newPDE}).

Analogously to the situation for ODEs, one has the following remarks.\\
\begin{rem}\label{rem:solsPDE}
If 
$(\alpha,\beta)=\left(A(x_1,x_2),B(x_1,x_2)\right)$ solves PDE system (\ref{reducnewPDE}), then due to the satisfaction of the integrability condition in PDE system (\ref{reducnewPDE}), $\frac{\partial \alpha}{\partial x_2}=\frac{\partial \beta}{\partial x_1}$, there exists a function $U(x_1,x_2)$ solving $\left( \tfrac{\partial U}{\partial x_1}(x_1,x_2), \tfrac{\partial U}{\partial x_2}(x_1,x_2) \right)=\left(A(x_1,x_2),B(x_1,x_2)\right)$. More importantly
$$u=U(x_1,x_2)+C$$
solves PDE (\ref{newPDE}) for any constant $C$. Conversely, if $u=U(x_1,x_2)$ solves PDE (\ref{newPDE}), then $(\alpha,\beta)=\left( \tfrac{\partial U}{\partial x_1}(x_1,x_2), \tfrac{\partial U}{\partial x_2}(x_1,x_2) \right)$ solves PDE system (\ref{reducnewPDE}). But the relationship between the solutions of PDE (\ref{newPDE}) and PDE system (\ref{reducnewPDE}) is not one-to-one.

In particular, for any solution of PDE system (\ref{reducnewPDE}), there is an infinite number of solutions of PDE (\ref{newPDE}). Conversely, for any solution of PDE (\ref{newPDE}), there is a unique solution of PDE system (\ref{reducnewPDE}). Thus PDE (\ref{newPDE}) and PDE system (\ref{reducnewPDE}) are nonlocally related.\\
\end{rem}

\begin{example}
\normalfont
As an example, consider the second-order PDE from \cite{blumanyang2013} given by
\begin{equation}\label{reacdiffeq}
\frac{\partial^2 u}{\partial x_1^2}-\frac{\partial u}{\partial x_2}+u \ln u=0.
\end{equation}
PDE (\ref{reacdiffeq}) admits the point symmetry 
\begin{equation}\label{symmreacdiffeq} \mathbb{X}=e^{x_2} u \dfrac{\partial}{\partial u}.
\end{equation}
In terms of canonical coordinates
\begin{equation}
\label{canonreacdiffvar} r_1=x_1, \quad r_2=x_2, \quad s=e^{-x_2}\ln u,
\end{equation}
PDE (\ref{reacdiffeq}) becomes
\begin{equation}\label{canonreacdiffeq}  \dfrac{\partial^2 s}{\partial r_1^2} -\dfrac{\partial s}{\partial r_2} +e^{r_2} \left(\dfrac{\partial s}{\partial r_1} \right)^2 =0.
\end{equation}
Without loss of generality, PDE (\ref{canonreacdiffeq}) can be re-written in the form
\begin{equation}\label{canonreacdiffeqorig}  \dfrac{\partial^2 u}{\partial x_1^2}-\dfrac{\partial u}{\partial x_2} +e^{x_2} \left(\dfrac{\partial u}{\partial x_1} \right)^2 =0.
\end{equation}
Let $\alpha = \dfrac{\partial u}{\partial x_1}$, $\beta=\dfrac{\partial u}{\partial x_2}$. Then PDE (\ref{canonreacdiffeqorig}) ``reduces'' to the nonlocally related first-order PDE system
\begin{equation}\label{reduccanonreacdiffeq}
\left\{\begin{array}{l}
 \dfrac{\partial \alpha}{\partial x_1}  -\beta+e^{x_2} \alpha^2 =0, \vspace*{0.15cm} \\
\dfrac{\partial \alpha}{\partial x_2}=\dfrac{\partial \beta}{\partial x_1}.
\end{array}\right.
\end{equation}
Consider a particular solution of PDE system (\ref{reduccanonreacdiffeq}) given by
\begin{equation}
\label{solpartreduccanonreacdiffeq}
\alpha=-\dfrac{1}{2}x_1 e^{-x_2}, \quad \beta=\dfrac{1}{4} \left(x_1^2-2 \right) e^{-x_2}.
\end{equation}
Correspondingly
\begin{equation}
\label{solpartreacdiffeq}
u=\dfrac{1}{4} \left(2-x_1^2 \right)e^{-x_2}+C 
\end{equation}
solves PDE (\ref{canonreacdiffeqorig}) with $C$ an arbitrary constant. 

Now consider a particular solution of PDE (\ref{canonreacdiffeqorig}) given by
\begin{equation}
\label{solpartcanonreacdiffeqorig}
u=x_1+ e^{x_2}.
\end{equation}
Correspondingly
\begin{equation}
\label{solpartreduccanonreacdiffeq2}
\alpha=1, \quad \beta=e^{x_2},
\end{equation}
solves PDE system (\ref{reduccanonreacdiffeq}). This illustrates that the relationship between the solutions of PDE (\ref{canonreacdiffeqorig}) and PDE system (\ref{reduccanonreacdiffeq}) is not one-to-one.\\
\end{example}

\begin{rem}\label{rem:symmPDE}
From the relationship between the solutions of PDE (\ref{newPDE}) and PDE system (\ref{reducnewPDE}), any symmetry of PDE (\ref{newPDE}) is a symmetry of PDE system (\ref{reducnewPDE}); and vice versa, any symmetry of PDE system (\ref{reducnewPDE}) is a symmetry of PDE (\ref{newPDE}).

A point symmetry of PDE (\ref{newPDE}) is not necessarily a point symmetry of PDE system (\ref{reducnewPDE}). Conversely, a point symmetry of PDE system (\ref{reducnewPDE}) is not necessarily a point symmetry of PDE (\ref{newPDE}).\\
\end{rem}

To illustrate the three different symmetry situations that could arise, we consider a second-order PDE admitting more than one point symmetry such that: (1) a second point symmetry of the given PDE yields a point symmetry of the ``reduced'' PDE system obtained from its invariance under the first point symmetry; (2) a second point symmetry of the given PDE yields a symmetry which is not a point symmetry of the ``reduced'' PDE system obtained from its invariance under the first point symmetry (this symmetry is a nonlocal symmetry of the ``reduced'' PDE system); (3) a point symmetry of the ``reduced'' PDE system obtained from the invariance under a point symmetry of the given PDE yields a symmetry which is not a point symmetry of the given PDE (this symmetry is a nonlocal symmetry of the given PDE).\\

\begin{example}\label{ex:PDEbothdirec}
\normalfont
Consider the second-order PDE given by
\begin{equation}\label{nonldiffeq}
\dfrac{\partial u}{\partial x_2}= \left( \frac{\partial u}{\partial x_1} \right)^{-\frac{4}{3}}\dfrac{\partial^2 u}{\partial x_1^2}.
\end{equation}

The admitted point symmetries for PDE (\ref{nonldiffeq}) appear in \cite{ovsiannikov1959,bluman1967}. PDE (\ref{nonldiffeq}) admits five point symmetries given by
\begin{equation}\label{symmnonldiffeq}
\begin{array}{c}
\displaystyle \mathbb{X}_1= \frac{\partial}{\partial u}, \quad \mathbb{X}_2=x_1 \frac{\partial}{\partial x_1}+2 x_2 \frac{\partial}{\partial x_2}+ u \frac{\partial}{\partial u}, \quad \mathbb{X}_3=\frac{\partial}{\partial x_1}, \vspace*{0.25cm}\\ \displaystyle \mathbb{X}_4=\frac{\partial}{\partial x_2}, \quad \mathbb{X}_5=2 x_1 \frac{\partial}{\partial x_1}-u \frac{\partial}{\partial u}.
\end{array}
\end{equation}
\medskip

\noindent Consider $\mathbb{X}_1$. From the invariance of PDE (\ref{nonldiffeq}) under the translation symmetry $\mathbb{X}_1$, we can directly introduce $\alpha=u_1=\frac{\partial u}{\partial x_1}$, $\beta=u_2=\frac{\partial u}{\partial x_2}$. Correspondingly, PDE (\ref{nonldiffeq}) ``reduces'' to the nonlocally related first-order PDE system
\begin{equation}
\label{reducnonldiffsys}
\left\{ \begin{array}{l} \beta = \alpha^{-\frac{4}{3}} \dfrac{\partial \alpha}{\partial x_1}, \vspace*{0.25cm}\\ 
\dfrac{\partial \alpha}{\partial x_2}=\dfrac{\partial \beta}{\partial x_1}.
\end{array} \right.
\end{equation}
Then
$$\mathbb{X}_2^{(1)}=\mathbb{X}_2- u_2 \frac{\partial}{\partial u_2}.$$
Hence the point symmetry $\mathbb{X}_2$ in (\ref{symmnonldiffeq}) of PDE (\ref{nonldiffeq}) becomes the point symmetry
\begin{equation}\label{X2X1symmreducnonldiffsys}
\mathbb{X}_2=x_1 \frac{\partial}{\partial x_1}+2 x_2 \frac{\partial}{\partial x_2}- \beta \dfrac{\partial}{\partial \beta}
\end{equation}
of PDE system (\ref{reducnonldiffsys}).\\

\noindent Now consider $\mathbb{X}_2$. In terms of canonical coordinates
\begin{equation}
\label{canonvarnonldiffeq} r_1=\dfrac{x_2}{x_1^2}, \quad r_2=\dfrac{u}{x_1}, \quad s=\log x_1,
\end{equation}
PDE (\ref{nonldiffeq}) becomes
\begin{equation}\label{canon2nonldiffeq} \begin{array}{c}
-\dfrac{\partial s}{\partial r_1}\left( \dfrac{\partial s}{\partial r_2}\right)^2+ \left( \dfrac{1+r_2 \dfrac{\partial s}{\partial r_2}+2 r_1 \dfrac{\partial s}{\partial r_1}}{\dfrac{\partial s}{\partial r_2}} \right)^{-\frac{4}{3}} \Bigg( -4 r_1 \left( 1+2 r_1 \dfrac{\partial s}{\partial r_1}\right) \dfrac{\partial s}{\partial r_2} \dfrac{\partial^2 s}{\partial r_1 \partial r_2}    \vspace*{0.15cm}\\ 
\qquad  +\left( 1+2 r_1 \dfrac{\partial s}{\partial r_1}\right)^2 \dfrac{\partial^2 s}{\partial r_2^2} + \left(\dfrac{\partial s}{\partial r_2} \right)^2 \left( -1+2 r_1 \dfrac{\partial s}{\partial r_1} +4 r_1^2  \dfrac{\partial^2 s}{\partial r_1^2} \right)  \Bigg)= 0.
\end{array}
\end{equation}
Let
\begin{equation}\label{auxvar2nonldiffeq}
\begin{array}{rcl}
\alpha &=&\dfrac{\partial s}{\partial r_1}= -\dfrac{x_1^2 u_2}{x_1 u_1+2 x_2 u_2-u}, \vspace*{0.15cm}\\
\beta&=&\dfrac{\partial s}{\partial r_2}=\dfrac{x_1}{x_1 u_1+2 x_2 u_2-u}.
\end{array} 
\end{equation}
Then PDE (\ref{nonldiffeq}) ``reduces'' to the nonlocally related first-order PDE system
\begin{equation}
\label{reduc2nonldiffsys}
\left\{ \begin{array}{l} \alpha \beta^2- \left( \dfrac{1+r_2 \beta+2 r_1 \alpha}{\beta} \right)^{-\frac{4}{3}} \Bigg( -4 r_1 \left( 1+2 r_1 \alpha\right) \beta   \dfrac{\partial \beta}{\partial r_1}  \vspace*{0.15cm}\\ 
\qquad   + \left( 1+2 r_1 \alpha\right)^2 \dfrac{\partial \beta}{\partial r_2}+ \beta^2 \left( -1+2 r_1 \alpha +4 r_1^2  \dfrac{\partial \alpha}{\partial r_1} \right)  \Bigg)= 0, \vspace*{0.25cm}\\ 
\displaystyle \frac{\partial \alpha}{\partial r_2}=\frac{\partial \beta}{\partial r_1}.
\end{array} \right.
\end{equation}
After applying
$$\mathbb{X}_1^{(1)}=\mathbb{X}_1=\frac{\partial}{\partial u},$$
to (\ref{canonvarnonldiffeq}) and (\ref{auxvar2nonldiffeq}), one obtains 
$$\mathbb{X}_1^{(1)}r_1=0, \qquad \mathbb{X}_1^{(1)}r_2=e^{-s}, \qquad \mathbb{X}_1^{(1)}\alpha= e^{-s}\alpha \beta, \qquad \mathbb{X}_1^{(1)}\beta= e^{-s} \beta^2.$$
Hence the point symmetry $\mathbb{X}_1$ in (\ref{symmnonldiffeq}) of PDE (\ref{nonldiffeq}) becomes the symmetry
\begin{equation}\label{X1X2symmreducnonldiffsys}
\mathbb{X}_1=e^{-s} \frac{\partial}{\partial r_2}+e^{-s}\alpha \beta \frac{\partial}{\partial \alpha}+ e^{-s} \beta^2  \dfrac{\partial}{\partial \beta},
\end{equation}
of PDE system (\ref{reduc2nonldiffsys}). 

But the symmetry (\ref{X1X2symmreducnonldiffsys}) of PDE system (\ref{reduc2nonldiffsys}) is not a point symmetry of PDE system (\ref{reduc2nonldiffsys}) since $\mathbb{X}_1$ depends explicitly on $s$. In particular $\mathbb{X}_1$ is a nonlocal symmetry of 
PDE system (\ref{reduc2nonldiffsys}).

Consider again PDE system (\ref{reducnonldiffsys}). The admitted point symmetries for PDE system (\ref{reducnonldiffsys}) appear in \cite{ovsiannikov1959}. In particular, PDE system (\ref{reducnonldiffsys}) admits the point symmetry 
\begin{equation}
\label{nonlocalsymmdiffeq}
\mathbb{Y}=x_1^2 \frac{\partial}{\partial x_1}-3 x_1 \alpha \frac{\partial}{\partial \alpha}-\left( 3\alpha^{-\frac{1}{3}}+x_1 \beta \right) \dfrac{\partial}{\partial \beta}.
\end{equation}
In comparison with the point symmetries (\ref{symmnonldiffeq}) of PDE (\ref{nonldiffeq}), it directly follows that the point symmetry $\mathbb{Y}$ yields a symmetry of PDE (\ref{nonldiffeq}) which is not a point symmetry of PDE (\ref{nonldiffeq}), since the infinitesimal components corresponding to the variables $(x_1,x_2)$ of the admitted point symmetries of PDE (\ref{nonldiffeq}) are not the same as those for $\mathbb{Y}$.\\
\end{example}

\begin{rem}\label{rem:furtherNLRPDE}
\noindent
If a point symmetry of PDE (\ref{newPDE}) is also a point symmetry of PDE system (\ref{reducnewPDE}), one can construct explicitly another nonlocally related system of PDEs for PDE system (\ref{reducnewPDE}) in addition to the scalar PDE (\ref{newPDE}), and a further nonlocally related system of PDEs for PDE (\ref{newPDE}) \cite{blumanrafamsml}.\\
\end{rem}

\begin{example}\label{ex:furtherNLRPDEsys}
\normalfont
Consider again the second-order PDE given by (\ref{nonldiffeq}) and the ``reduced'' first-order PDE system (\ref{reducnonldiffsys}). Consider $\mathbb{X}_2$ given by (\ref{X2X1symmreducnonldiffsys}). From Example \ref{ex:PDEbothdirec}, we know that $\mathbb{X}_2$ is a point symmetry of PDE system (\ref{reducnonldiffsys}). In terms of canonical coordinates
\begin{equation}
\label{furthercanonvarnonldiffsys} r_1=x_1 \beta, \qquad  r_2=\dfrac{x_1^2}{x_2}, \qquad s_1=\log x_1, \qquad s_2=\alpha,
\end{equation}
PDE system (\ref{reducnonldiffsys}) becomes
\begin{equation}\label{furthercanonnonldiffsys}
\left\{ \begin{array}{l}
r_2^2 \left( \dfrac{\partial s_1}{\partial r_2} \dfrac{\partial s_2}{\partial r_1}-\dfrac{\partial s_1}{\partial r_1} \dfrac{\partial s_2}{\partial r_2} \right)+r_1 \dfrac{\partial s_1}{\partial r_1}+2 r_2 \dfrac{\partial s_1}{\partial r_2}-1=0, \vspace*{0.25cm}\\
2 r_2 \left( \dfrac{\partial s_1}{\partial r_2} \dfrac{\partial s_2}{\partial r_1}-\dfrac{\partial s_1}{\partial r_1} \dfrac{\partial s_2}{\partial r_2} \right)-\dfrac{\partial s_2}{\partial r_1}+r_1  s_2^{\frac{4}{3}}\dfrac{\partial s_1}{\partial r_1} =0.
\end{array} \right.
\end{equation}
Let $\delta=\dfrac{\partial s_1}{\partial r_1}$, $\omega=\dfrac{\partial s_1}{\partial r_2}$. Then PDE system (\ref{furthercanonnonldiffsys}) becomes the ``reduction of order'' PDE system
\begin{equation}\label{further2canonnonldiffsys}
\left\{ \begin{array}{l}
r_2^2 \left( \omega \dfrac{\partial s_2}{\partial r_1}-\delta \dfrac{\partial s_2}{\partial r_2} \right)+r_1 \delta+2 r_2 \omega-1=0, \vspace*{0.25cm}\\
2 r_2 \left( \omega \dfrac{\partial s_2}{\partial r_1}-\delta \dfrac{\partial s_2}{\partial r_2} \right)-\dfrac{\partial s_2}{\partial r_1} +r_1  s_2^{\frac{4}{3}}\delta=0,
\vspace*{0.25cm}\\
\displaystyle \frac{\partial \delta}{\partial r_2}=\dfrac{\partial \omega}{\partial r_1},
\end{array} \right.
\end{equation}
nonlocally related to PDE system (\ref{reducnonldiffsys}) and PDE (\ref{nonldiffeq}).
\end{example}

\subsection{A PDE with three independent variables}\label{pde3}
Consider an $n$th-order PDE with three independent variables $\underline{x}=(x_1,x_2,x_3)$ and dependent variable $u$ given by
\begin{equation}
\label{PDE3var}
K\left(\underline{x},u,\partial u,\ldots,\partial^n u \right)=0,
\end{equation}
with $\partial^l u$ denoting its partial derivatives of $u$ of order $l$, $l=1,2,\ldots,n$.

Assume that PDE (\ref{PDE3var}) admits a point symmetry with infinitesimal generator
\begin{equation}
\label{symmPDE3var}
\mathbb{X}=\xi_1 \left(\underline{x},u\right)\dfrac{\partial}{\partial x_1}+\xi_2 \left(\underline{x},u\right)\dfrac{\partial}{\partial x_2}+\xi_3 \left(\underline{x},u\right)\dfrac{\partial}{\partial x_3}+\eta\left(\underline{x},u\right)\dfrac{\partial}{\partial u}
\end{equation}
whose $n$th extension is given by
\begin{equation}
\label{symmPDEext3var}
\mathbb{X}^{(n)}=\mathbb{X}+ \eta_j^{(1)}\left(\underline{x},u,\partial u\right)\dfrac{\partial}{\partial u_j}+\cdots+\eta_{j_1 j_2 \cdots j_n}^{(n)}\left(\underline{x},u,\partial u,\ldots,\partial^n u \right)\dfrac{\partial}{\partial u_{j_1 j_2 \cdots j_n}},
\end{equation}
where $u_j=\frac{\partial u}{\partial x_j}$, $j=1,2,3$, and $j_k=1,2,3,$ for $k=1,\ldots,n$, and the $k$th-extended infinitesimals  $\eta^{(k)}$ satisfy the recursion relations
$$\begin{array}{rcl}\eta_j^{(1)}&=& D_j \eta- \left( D_j \xi_i \right)u_i, \quad j=1,2,3, \vspace*{0.05cm}\\ 
\eta_{j_1 j_2 \cdots j_k}^{(k)}&=&  D_{j_k} \eta_{j_1 j_2 \cdots j_{k-1}}^{(k-1)}- \left( D_{j_k} \xi_i \right)u_{j_1 j_2 \cdots j_{k-1} i},
\end{array}$$
in terms of the total derivative operators given by (\ref{totalderivative}).

Then $\mathbb{X}$ is a point symmetry of PDE (\ref{PDE3var}) if and only if 
$$ \mathbb{X}^{(n)}K\left(\underline{x},u,\partial u,\ldots,\partial^n u \right) =0 \quad \mbox{when} \quad K\left(\underline{x},u,\partial u,\ldots,\partial^n u \right)=0.$$
We introduce canonical coordinates, $\underline{r}(\underline{x},u)=\left(r_1(\underline{x},u),r_2(\underline{x},u),r_3(\underline{x},u)\right)$ and $s(\underline{x},u)$, for $\mathbb{X}$, satisfying
$$\mathbb{X} \underline{r}=0, \quad  \mathbb{X}s=1.$$
Thus, equivalently, PDE (\ref{PDE3var}) is invariant under translations in $s$, i.e., PDE (\ref{PDE3var}) is invariant under the point symmetry (\ref{symmPDE3var}) now written as $\mathbb{X}=\dfrac{\partial}{\partial s}$.\vspace*{0.1cm}

In terms of corresponding canonical coordinates, PDE (\ref{PDE3var}) invertibly transforms by a point transformation to a PDE of the form
\begin{equation}
\label{transformPDE3var}
\widehat{K}\left(\underline{r},\partial s,\ldots,\partial^n s \right)=0.
\end{equation}

Hence, without loss of generality, if a given PDE with three independent variables admits a point symmetry, after using corresponding canonical coordinates, it can be re-written in the form
\begin{equation}
\label{newPDE3var}
L\left(\underline{x},\partial u,\ldots,\partial^n u \right)=0,
\end{equation}
i.e., PDE (\ref{newPDE3var}) is invariant under the translation symmetry $\mathbb{X}=\dfrac{\partial}{\partial u}$.\vspace*{0.1cm}

From the form of PDE (\ref{newPDE3var}), we can naturally introduce as auxiliary dependent variables $\underline{\alpha}= \nabla u=\left(\frac{\partial u}{\partial x_1},\frac{\partial u}{\partial x_2},\frac{\partial u}{\partial x_3}\right)$. Then PDE (\ref{newPDE3var}) becomes the ``reduction of order'' PDE system
\begin{equation}
\label{reducnewPDE3var}
\left\{ \begin{array}{l} L\left(\underline{x},\underline{\alpha}, \partial \underline{\alpha}, \ldots,\partial^{n-1} \underline{\alpha} \right) =0, \vspace*{0.15cm}\\ 
\mbox{curl} \, \underline{\alpha}= \nabla \times \underline{\alpha}=0,
\end{array} \right.
\end{equation}
nonlocally related to PDE (\ref{newPDE3var}).

Analogously to the situation for PDEs with two independent variables, one has the following remarks.\\

\begin{rem}\label{rem:solsPDE3var}
If 
$\underline{\alpha}=\underline{A}\left(\underline{x}\right)$ solves PDE system (\ref{reducnewPDE3var}), then due to the satisfaction of the integrability condition in PDE system (\ref{reducnewPDE3var}), $\nabla \times \underline{\alpha}=0$, there exists a function $U(\underline{x})$ solving $\nabla U (\underline{x})=\underline{A}\left(\underline{x}\right)$. More importantly
$$u=U(\underline{x})+C$$
solves PDE (\ref{newPDE3var}) for any constant $C$. Conversely, if $u=U(\underline{x})$ solves PDE (\ref{newPDE3var}), then $\underline{\alpha}=\nabla U (\underline{x})$ solves PDE system (\ref{reducnewPDE3var}). But the relationship between the solutions of PDE (\ref{newPDE3var}) and PDE system (\ref{reducnewPDE3var}) is not one-to-one.

In particular, for any solution of PDE system (\ref{reducnewPDE3var}), there is an infinite number of solutions of PDE (\ref{newPDE3var}). Conversely, for any solution of PDE (\ref{newPDE3var}), there is a unique solution of PDE system (\ref{reducnewPDE3var}). Thus PDE (\ref{newPDE3var}) and PDE system (\ref{reducnewPDE3var}) are nonlocally related.\\
\end{rem}

\begin{rem}\label{rem:symmPDE3var}
From the relationship between the solutions of PDE (\ref{newPDE3var}) and PDE system (\ref{reducnewPDE3var}), any symmetry of PDE (\ref{newPDE3var}) is a symmetry of PDE system (\ref{reducnewPDE3var}); and vice versa, any symmetry of PDE system (\ref{reducnewPDE3var}) is a symmetry of PDE (\ref{newPDE3var}).

A point symmetry of PDE (\ref{newPDE3var}) is not necessarily a point symmetry of PDE system (\ref{reducnewPDE3var}). Conversely, a point symmetry of PDE system (\ref{reducnewPDE3var}) is not necessarily a point symmetry of PDE (\ref{newPDE3var}).\\
\end{rem}

\begin{rem}\label{rem:furtherNLRPDE3var}
\noindent
If a point symmetry of PDE (\ref{newPDE3var}) is also a point symmetry of PDE system (\ref{reducnewPDE3var}), one can construct explicitly another nonlocally related system of PDEs for PDE system (\ref{reducnewPDE3var}) in addition to the scalar PDE (\ref{newPDE3var}), and a further nonlocally related system of PDEs for PDE (\ref{newPDE3var}) \cite{blumanrafamsml,wang2023}.
\end{rem}

\subsection{A PDE with more than three independent variables}\label{pdep}
Consider an $n$th-order PDE with independent variables $\underline{x}=\left(x_1,x_2,\ldots,x_p\right)$, with $p>3$, and dependent variable $u$ given by
\begin{equation}
\label{PDEpvar}
K\left(\underline{x},u,\partial u,\ldots,\partial^n u \right)=0,
\end{equation}
with $\partial^l u$ denoting its partial derivatives of $u$ of order $l$, $l=1,2,\ldots,n$.

Assume that PDE (\ref{PDEpvar}) admits a point symmetry with infinitesimal generator
\begin{equation}
\label{symmPDEpvar}
\mathbb{X}=\xi_j \left(\underline{x},u\right)\dfrac{\partial}{\partial x_j}+\eta\left(\underline{x},u\right)\dfrac{\partial}{\partial u}
\end{equation}
whose $n$th extension is given by
\begin{equation}
\label{symmPDEextpvar}
\mathbb{X}^{(n)}=\mathbb{X}+ \eta_j^{(1)}\left(\underline{x},u,\partial u\right)\dfrac{\partial}{\partial u_j}+\cdots+\eta_{j_1 j_2 \cdots j_n}^{(n)}\left(\underline{x},u,\partial u,\ldots,\partial^n u \right)\dfrac{\partial}{\partial u_{j_1 j_2 \cdots j_n}},
\end{equation}
where $u_j=\frac{\partial u}{\partial x_j}$, $j=1,\ldots,p$, and $j_k=1,\ldots,p$, for $k=1,\ldots,n$, and the $k$th-extended infinitesimals  $\eta^{(k)}$ satisfy the recursion relations
$$\begin{array}{rcl}\eta_j^{(1)}&=& D_j \eta- \left( D_j \xi_i \right)u_i, \quad j=1,\ldots,p, \vspace*{0.05cm}\\ 
\eta_{j_1 j_2 \cdots j_k}^{(k)}&=&  D_{j_k} \eta_{j_1 j_2 \cdots j_{k-1}}^{(k-1)}- \left( D_{j_k} \xi_i \right)u_{j_1 j_2 \cdots j_{k-1} i},
\end{array}$$
in terms of the total derivative operators given by (\ref{totalderivative}).

Then $\mathbb{X}$ is a point symmetry of PDE (\ref{PDEpvar}) if and only if 
$$ \mathbb{X}^{(n)}K\left(\underline{x},u,\partial u,\ldots,\partial^n u \right) =0 \quad \mbox{when} \quad K\left(\underline{x},u,\partial u,\ldots,\partial^n u \right)=0.$$
We introduce canonical coordinates, $\underline{r}(\underline{x},u)=\left(r_1(\underline{x},u),\ldots,r_p(\underline{x},u)\right)$ and $s(\underline{x},u)$, for $\mathbb{X}$, satisfying
$$\mathbb{X}\underline{r}=0, \quad  \mathbb{X}s=1.$$
Thus, equivalently, PDE (\ref{PDEpvar}) is invariant under translations in $s$, i.e., PDE (\ref{PDEpvar}) is invariant under the point symmetry (\ref{symmPDEpvar}) now written as $\mathbb{X}=\dfrac{\partial}{\partial s}$.\vspace*{0.1cm}

In terms of corresponding canonical coordinates, PDE (\ref{PDEpvar}) invertibly transforms by a point transformation to a PDE of the form
\begin{equation}
\label{transformPDEpvar}
\widehat{K}\left(\underline{r},\partial s,\ldots,\partial^n s \right)=0.
\end{equation}

Hence, without loss of generality, if a given PDE with $p$ independent variables, $p>3$, admits a point symmetry, after using corresponding canonical coordinates, it can be re-written in the form
\begin{equation}
\label{newPDEpvar}
L\left(\underline{x},\partial u,\ldots,\partial^n u \right)=0,
\end{equation}
i.e., PDE (\ref{newPDEpvar}) is invariant under the translation symmetry $\mathbb{X}=\dfrac{\partial}{\partial u}$.\vspace*{0.1cm}

From the form of PDE (\ref{newPDEpvar}), we can naturally introduce as auxiliary dependent variables 
$$\underline{\alpha}=\nabla u=\left(\frac{\partial u}{\partial x_1},\ldots,\frac{\partial u}{\partial x_p}\right).$$
Then PDE (\ref{newPDEpvar}) becomes the ``reduction of order'' PDE system
\begin{equation}
\label{reducnewPDEpvar}
\left\{ \begin{array}{l} L\left(\underline{x},\underline{\alpha}, \partial \underline{\alpha}, \ldots,\partial^{n-1} \underline{\alpha} \right) =0, \vspace*{0.15cm}\\ 
\dfrac{\partial \alpha_i}{\partial x_j}- \dfrac{\partial \alpha_j}{\partial x_i}=0, \quad i=1,2,\ldots,p-1, \quad j=i+1,\ldots, p,
\end{array} \right.
\end{equation}
nonlocally related to PDE (\ref{newPDEpvar}).

Analogously to the situations for PDEs with two and three independent variables, one has the following remarks.\\

\begin{rem}\label{rem:solsPDEpvar} 
If 
$\underline{\alpha}=\underline{A}\left(\underline{x}\right)$ solves PDE system (\ref{reducnewPDEpvar}), then due to the satisfaction of the integrability conditions in PDE system (\ref{reducnewPDEpvar}), there exists a function $U(\underline{x})$ solving $\nabla U (\underline{x})=\underline{A}\left(\underline{x}\right)$. More importantly
$$u=U(\underline{x})+C$$
solves PDE (\ref{newPDEpvar}) for any constant $C$. Conversely, if $u=U(\underline{x})$ solves PDE (\ref{newPDEpvar}), then $\underline{\alpha}=\nabla U (\underline{x})$ solves PDE system (\ref{reducnewPDEpvar}). But the relationship between the solutions of PDE (\ref{newPDEpvar}) and PDE system (\ref{reducnewPDEpvar}) is not one-to-one.

In particular, for any solution of PDE system (\ref{reducnewPDEpvar}), there is an infinite number of solutions of PDE (\ref{newPDEpvar}). Conversely, for any solution of PDE (\ref{newPDEpvar}), there is a unique solution of PDE system (\ref{reducnewPDEpvar}). Thus PDE (\ref{newPDEpvar}) and PDE system (\ref{reducnewPDEpvar}) are nonlocally related \cite{blumanyang2013,yang2013,blumanyang2014,DIpaper1}.\\
\end{rem}

\begin{rem}\label{rem:symmPDEpvar}
From the relationship between the solutions of PDE (\ref{newPDEpvar}) and PDE system (\ref{reducnewPDEpvar}), any symmetry of PDE (\ref{newPDEpvar}) is a symmetry of PDE system (\ref{reducnewPDEpvar}); and vice versa, any symmetry of PDE system (\ref{reducnewPDEpvar}) is a symmetry of PDE (\ref{newPDEpvar}).

A point symmetry of PDE (\ref{newPDEpvar}) is not necessarily a point symmetry of PDE system (\ref{reducnewPDEpvar}). Conversely, a point symmetry of PDE system (\ref{reducnewPDEpvar}) is not necessarily a point symmetry of PDE (\ref{newPDEpvar}).\\
\end{rem}

\begin{rem}\label{rem:furtherNLRPDEpvar}
\noindent
If a point symmetry of PDE (\ref{newPDEpvar}) is also a point symmetry of PDE system (\ref{reducnewPDEpvar}), one can construct explicitly another nonlocally related system of PDEs for PDE system (\ref{reducnewPDEpvar}) in addition to the scalar PDE (\ref{newPDEpvar}), and a further nonlocally related system of PDEs for PDE (\ref{newPDEpvar}) \cite{blumanrafamsml}.
\end{rem}

\section{Concluding Remarks}\label{sec:concluding}

The following concluding remarks are noted.

\subsection{The extension to systems of DEs with more than one dependent variable}

For simplicity, the theory presented in Sections \ref{sec:liereducODE} and \ref{sec:PDE} focused on scalar ODEs and scalar PDEs. Nevertheless, the concepts and situations described in these sections are easily extended to systems of ODEs and systems of PDEs. In particular, in Example \ref{ex:furtherNLRPDEsys} we applied the symmetry-based method
to the nonlocally related PDE system (\ref{reducnonldiffsys}) of PDE (\ref{nonldiffeq}). For further details on the use of the symmetry-based method to construct nonlocally related systems for a DE system, with two or more dependent variables,
see \cite{blumanyang2013,yang2013,blumanyang2014,DIpaper1}.

\subsection{``Reduction of Order''}

As we have shown, if a PDE system admits a point symmetry, one can construct a nonlocally related ``reduced'' system of PDEs for the given PDE system. In this case, the ``reduction of order'' is a reduction of order in one of its dependent variables, namely, a reduction of order with respect to the translated dependent variable that arises from the canonical coordinates of the point symmetry.

In this paper we have shown that the symmetry-based method is a natural extension of Lie's reduction of order algorithm for ODEs.  Hence it is more appropriate to replace the phrase ``nonlocally related inverse potential system'' appearing in the preceding references \cite{blumanyang2013,yang2013,blumanyang2014,blumanrafamsml,blumanzuhal,DIpaper1,DIpaper2,wang2023} by the phrase ``nonlocally related reduced PDE system''.  Moreover, there is no longer a need to consider a ``locally related intermediate system'' that appeared in the same preceding references.

\subsection{Reduction of order algorithm for DEs admitting a solvable Lie group of point transformations: the natural extension from ODEs to PDEs}

In Section \ref{sec:liereducODE}, we reviewed Lie's reduction of order algorithm showing that the order of an $n$th-order ODE admitting a point symmetry can be reduced constructively by one through the use of canonical coordinates. Consequently,  an infinite number of solutions of a given $n$th-order ODE is obtained from a solution of a nonlocally related reduced $(n-1)$th-order ODE through quadrature. If an $n$th-order ODE, $n\geq2$, admits precisely two point symmetries, then the order can be reduced constructively by two if the point symmetries are used in the correct order. But if an $n$th-order ODE admits $q$ point symmetries, $n\geq2$ and $q>2$, it does not necessarily follow that the order can be reduced by more than one. However, if an $n$th-order ODE admits $q$ point symmetries whose infinitesimal generators form a $q$-dimensional solvable Lie algebra, then there exists a solvable chain of $q$ subalgebras whose successive subalgebras contain $1,\ldots,q$ point symmetries. Consequently, a chain of nonlocally related reduced order ODEs of orders $n-1,n-2,\ldots,n-q$, assuming $q \leq n$, is obtained if the point symmetries are used in the correct order. Moreover, from a solution of any such  nonlocally related reduced order ODE, one obtains an infinite number of solutions of the given ODE through quadrature(s). The use of solvable groups for successive reduction of order appeared in \cite{bianchi1918,eisenhart1933,olver1986,blumankumei,bluman1990}.

As examples, consider ODEs (\ref{secondex}) and (\ref{canonblasiusorig}).  ODE (\ref{secondex}) admits two point symmetries (\ref{symmsecondex}) whose infinitesimal generators yield a two-dimensional Lie algebra with the commutation relation $\left[\mathbb{X}_1,\mathbb{X}_2\right]=-2\mathbb{X}_1$, whereas ODE (\ref{canonblasiusorig}) admits two point symmetries (\ref{symmblasiusorig}) whose infinitesimal generators similarly yield a two-dimensional Lie algebra with the commutation relation $\left[\mathbb{X}_1,\mathbb{X}_2\right]=-\mathbb{X}_1$. These are solvable algebras with commutation relations that immediately indicate that the nonlocally related reduced ODE arising from the admitted point symmetry $\mathbb{X}_1$ inherits the symmetry $\mathbb{X}_2$ as a point symmetry so that a further nonlocally related ODE arises. On the other hand, these commutation relations indicate that although the nonlocally related reduced ODE arising from the point symmetry $\mathbb{X}_2$ inherits the symmetry $\mathbb{X}_1$, it will not inherit $\mathbb{X}_1$ as a point symmetry and hence no further nonlocally related ODE can arise. These situations are shown explicitly for ODEs (\ref{secondex}) and (\ref{canonblasiusorig}). Further examples involving higher dimensional solvable algebras appear in \cite{bluman1990,blumankumei1989}.

In Section \ref{sec:PDE}, we naturally extended Lie's reduction of order algorithm yielding a nonlocally related ODE for an ODE admitting a point symmetry to a reduction of order algorithm yielding a nonlocally related PDE system for a PDE admitting a point symmetry. Correspondingly, if a given PDE system admits $q$ point symmetries whose infinitesimal generators form a $q$-dimensional solvable Lie algebra then one can naturally construct further nonlocally related reduced order PDE systems by using the point symmetries in a correct order from the resulting solvable chain of $q$ subalgebras.  For details, see \cite{blumanrafamsml}.

As an example, consider PDE (\ref{nonldiffeq}). PDE (\ref{nonldiffeq}) admits five point symmetries (\ref{symmnonldiffeq}) whose infinitesimal  generators $\mathbb{X}_1$ and $\mathbb{X}_2$ yield a two-dimensional subalgebra with the commutation relation $\left[\mathbb{X}_1,\mathbb{X}_2\right]=\mathbb{X}_1$. This immediately indicates that the nonlocally related PDE system arising from the admitted point symmetry $\mathbb{X}_1$ inherits the symmetry $\mathbb{X}_2$ as a point symmetry so that a further nonlocally related PDE system arises. On the other hand, this commutation relation indicates that although the nonlocally related PDE system arising from the point symmetry $\mathbb{X}_2$ inherits  the symmetry $\mathbb{X}_1$, it will not inherit $\mathbb{X}_1$ as a point symmetry and hence no further nonlocally related PDE system can arise.  These situations are shown explicitly for PDE (\ref{nonldiffeq}). Further examples involving higher dimensional solvable algebras appear in \cite{blumanrafamsml,wang2023}.\\

\subsection{Use of Nonlocally Related Systems}

Usually, a good mathematical method that is useful for a given DE, e.g., a symmetry method, numerical method, perturbation method, or qualitative method, still works subject to any invertible point transformation of its variables. But if such a method is not useful for a given DE, then it might be useful when applying the method to an equivalent nonlocally related DE since the nonlocally related DE is not equivalent to the given DE through any invertible point transformation. In this paper for many given DEs we have given examples of obtaining further concrete symmetries (nonlocal symmetries) that arise from  point symmetries of nonlocally related reduced systems.
 
\section*{Acknowledgments}
We thank the Natural Sciences and Engineering Research Council of Canada for financial support. G.W.B. thanks Subhankar Sil for a discussion that led to the writing of this paper.

\end{document}